\documentclass[11pt,reqno]{article}

\usepackage[usenames]{color}
\usepackage{amssymb}
\usepackage{amsmath}
\usepackage{amsthm}
\usepackage{amsfonts}
\usepackage{amscd}
\usepackage{graphicx}

\usepackage[colorlinks=true,
linkcolor=webgreen,
filecolor=webbrown,
citecolor=webgreen]{hyperref}

\definecolor{webgreen}{rgb}{0,.5,0}
\definecolor{webbrown}{rgb}{.6,0,0}

\usepackage{color}
\usepackage{fullpage}
\usepackage{float}

\usepackage{psfig}
\usepackage{graphics}
\usepackage{latexsym}
\usepackage{epsf}
\usepackage{breakurl}
\usepackage{enumitem}

\newcommand{\Li }{\operatorname{Li}}
\newcommand{\diff }{\,\mathrm{d}}

\begin{document}
	\sloppy
	\begin{center}
		\epsfxsize=4in
	\end{center}

	\newtheorem{theorem}{Theorem}           
	\newtheorem{lemma}{Lemma}               
	\newtheorem{corollary}{Corollary}
	\newtheorem{remark}{Remark}
	\newtheorem{example}{Example}

	\begin{center}
	\vskip 1cm{\LARGE \bf Explicit Evaluations of Euler Sums Involving Harmonic Numbers with Rational Arguments 
	}
	\vskip 1cm {\large
		Ali Olaikhan\\
		Grand Canyon University\\
		Phoenix, Arizona, United States\\
		\href{mailto:alishathri@yahoo.com}{\tt alishathri@yahoo.com}}
\end{center}

\vskip .2in
\begin{abstract}
	This paper derives closed-form expressions for four-parameter Euler sums containing generalized harmonic numbers with rational arguments, expressed in terms of the polylogarithm and Lerch transcendent, and reducible (under admissible parameter choices) to Riemann zeta and Hurwitz zeta values. It further obtains closed forms for related polylogarithmic integrals.

\end{abstract}

\section{Introduction and preliminaries}

In \cite[Corollary~2.4]{sofo2019}, Sofo established the remarkable identity
\begin{align}
	\sum_{k=1}^\infty\frac{H_{k/n}^{(p)}}{k^q}
	=&\,
	{(-1)}^{q}
	\sum_{r=1}^{p}
	\frac{1}{n^{q-r}}
	\binom{p+q-r-1}{p-r}
	\sum_{k=1}^{\infty}
	\frac{H_{n k}^{(r)}}{k^{p+q-r}}\notag\\
	\
	&
	-{(-1)}^{q}
	\sum_{r=2}^{p}
	\frac{1}{n^{q-r}}
	\binom{p+q-r-1}{p-r}
	\zeta(r)\zeta(p+q-r)\notag\\
	\
	&
	-{(-1)}^{q}
	\sum_{\ell=2}^{q-1}
	\frac{{(-1)}^{\ell}}{n^{q-\ell}}
	\binom{p+q-\ell-1}{q-\ell}
	\zeta(\ell)\zeta(p+q-\ell),
	\label{sofo}
\end{align}
where $p,n\in\mathbb{Z}_{\geqslant 1}$, $q\in\mathbb{Z}_{\geqslant 2}$, and $p+q$ is odd. Here, $H_n^{(p)}$ denotes the generalized harmonic numbers of order $p$, defined by
\[
H_n^{(p)}=
\sum_{k=1}^{n}\frac{1}{k^p},\quad p,n\in\mathbb{Z}_{\geqslant 1},
\]
and $\zeta(s)$ is the Riemann zeta function \cite[\S 25.2.1]{hurwitz}, defined by 
\[
\zeta(s)=
\sum_{n=1}^{\infty}\frac{1}{n^s},
\qquad
\Re(s)>1.
\]

Sofo also obtained an explicit evaluation of \eqref{sofo} in the particular case n=2, which was later derived independently in a simpler closed form in \cite[Theorem~6.1]{olaikhan}. Recently, the Euler sum $\sum_{k=1}^{\infty} H_{nk}^{(r)}/k^{p+q-r}$ appearing in \eqref{sofo} was evaluated in closed form \cite[Corollary~1]{olaikhan2026}, thereby yielding an explicit evaluation of $\sum_{k=1}^{\infty}H_{k/n}^{(p)}/k^q$.

Motivated by identity \eqref{sofo}, we derive in this paper closed-form expressions for the more general Euler sum (see Theorem \ref{thm})
\[
\sum_{k=1}^{\infty}
\frac{a^kH_{k/n}^{(p)}}{k^q}.
\qquad
a=\pm1,
\]
and the related logarithmic integral (see Theorem~\ref{thm2})
\begin{equation*}
	\int_0^1\frac{\ln^{p-1}(x)\Li_q(ax)}{x(1-x^n)}\,\diff x,
\end{equation*}
where $\Li_p(s)$ is the polylogarithm function, defined by \cite[p.~189, Eq.~(7.1)]{lewin}	
\[
\Li_p(s)=
\sum_{n=1}^{\infty}\frac{s^n}{n^p},
\qquad
p\in\mathbb{C}
\ \text{when}\ |s|<1,
\quad
\Re(p)>1
\ \text{when}\ |s|=1.
\]

These expressions are given in terms of the polylogarithm and Lerch transcendent, defined by \cite[p.~194, Eq.~(1)]{choibook}
\[
\Phi(x,s,t)=
\sum_{n=0}^{\infty}
\frac{x^n}{{(n+t)}^s},
\]
\[
t\notin\mathbb{Z}_{\leqslant0},
\qquad
s\in\mathbb{R}\ \text{when}\ |x|<1,
\qquad
s>1\ \text{when}\ |x|=1,
\]
and reducible (under admissible parameter choices) to Riemann zeta and Hurwitz zeta values, as illustrated in Section~\ref{ex}. The Hurwitz zeta function is defined by \cite[p.~15, Eq.~(1.3.1)]{andrews} 
\begin{equation*}
	\zeta(s,t) = \sum_{n=0}^\infty \frac{1}{{(n+t)}^s}, \quad \Re(s) > 1, \; t\notin\mathbb{Z}_{\leqslant 0}.
\end{equation*}

\section{Lemmas} 

\begin{lemma}\label{partial}
	Let $p,q,n\in\mathbb{Z}_{\geqslant 1}$. Then
	\begin{equation*}
		\frac{1}{k^q{(k+nj)}^{p}}
		=
		\sum_{r=1}^{p}
		\frac{{(-1)}^{q}
			\displaystyle\binom{p+q-r-1}{p-r}}{(nj)^{p+q-r}{(k+nj)}^r}+
		\sum_{\ell=1}^{q}
		\frac{ {(-1)}^{\ell+q} \displaystyle\binom{p+q-\ell-1}{q-\ell}
		}{(nj)^{p+q-\ell}k^\ell}.
	\end{equation*}
\end{lemma}
\begin{proof}
	Set $x=-k/nj$ in the partial fraction decomposition identity \cite[Eq.~(2.4)]{xu}
	\[
	\frac{1}{x^q{(1-x)}^p}
	=
	\sum_{r=1}^{p}
	\frac{\displaystyle\binom{p+q-r-1}{p-r}}{(1-x)^r}+
	\sum_{\ell=1}^{q}
	\frac{\displaystyle\binom{p+q-\ell-1}{q-\ell}}{x^\ell}.
	\]
\end{proof}

\begin{lemma} \label{m} Let $p,n\in\mathbb{Z}_{\geqslant 1}$ and $a=\pm 1$, $p\ne 1$ when $a=1$. Then
	\begin{equation*}
		\sum_{k=1}^\infty\frac{a^k}{{(k+n)}^p}=a^n\Li_p(a)-a^{n+1}\mathbf{H}_{n,a}^{(p)},
	\end{equation*}
	where $\Li_p(\cdot)$ is the polylogarithm and $\mathbf{H}_{n,a}^{(p)}$ is given by
	\begin{equation}
		\mathbf{H}_{n,a}^{(p)}=\sum_{j=1}^{n}\frac{a^{j-1}}{j^p},\quad \mathbf{H}_{n,a}^{(1)}:=\mathbf{H}_{n,a}, \label{def of H_n}
	\end{equation}
	where $a=1$ corresponds to the generalized harmonic numbers of order $p$, denoted by $H_n^{(p)}$, and $a=-1$ corresponds to the generalized skew-harmonic numbers of order $p$, denoted by $\overline{H}_n^{(p)}$.
\end{lemma}
\begin{proof} Setting the index $k=j-n$, we have
	\begin{align*}
		\sum_{k=1}^\infty\frac{a^k}{{(k+n)}^p}&=a^n\sum_{j=n+1}^\infty\frac{a^j}{j^p}\\
		&=a^n\left(\sum_{j=1}^\infty\frac{a^j}{j^p}-\sum_{j=1}^n\frac{a^j}{j^p}\right)\\
		&=a^n\left(\Li_p(a)-a\mathbf{H}_{n,a}^{(p)}\right).
	\end{align*}
\end{proof}

\begin{lemma}\label{mine}
	Let $p,n\in\mathbb{Z}_{\geqslant 1}$, $q\in\mathbb{Z}_{\geqslant 2}$, $a=\pm 1$, $b=\pm 1$, and $p+q$ odd. Then
	\begin{align*}
		\sum_{k=1}^\infty \frac{a^k \,  \mathbf{H}_{nk,b}^{(p)}}{k^q}=&\,\frac{b}{2} \left(1 + {(-1)}^q \right) \Li_q(a) \Li_p(b) + \frac{b}{2 n^p} \Li_{p+q}(a b^{n})  \\
		& -{(-1)}^{q} n^q b \sum_{j=0}^{\lfloor \frac{q}{2} \rfloor} \binom{p+q-2j-1}{p-1} 
		{n}^{-2j}\Li_{2j}(a b^{n}) \Li_{p+q-2j}(b)  \\
		& + \frac{{(-1)}^pb}{n^p} \sum_{j=0}^{\lfloor \frac{p}{2} \rfloor} \binom{p+q-2j-1}{q-1} 
		\Li_{2j}(a b^{n}) \Li_{p+q-2j}(a)  \\
		& + \frac{{(-1)}^pa b}{2 n^p} \sum_{j=1}^{n-1} \sum_{k=1}^p \binom{p+q-k-1}{q-1} 
		b^{j} {(-1)}^k\, \Theta(a b^{n}, k, j, n)\\
		&\quad\times\Phi\!\left(a, p+q-k, \frac{n-j}{n}\right),
	\end{align*}
	where $\mathbf{H}_{nk,b}^{(p)}$ is defined in \eqref{def of H_n}, $\lfloor\cdot\rfloor$ is the floor function, $\Li_p(\cdot)$ is the polylogarithm, $\Phi$ is the Lerch transcendent, and $\Theta$ is given by
	\begin{equation*}
		\Theta(c,q,s,r)=\Phi\!\left(c,q,\frac{s}{r}\right)+c{(-1)}^{q} \,\Phi\!\left(c,q,\frac{r-s}{r}\right),\quad \Theta(1,1,s,r) = \pi \cot\!\left(\frac{\pi s}{r}\right).
	\end{equation*} 
\end{lemma}

\begin{proof}
	Given in \cite[Corollary~1]{olaikhan2026}.
\end{proof}

The case $(p,q,n,a,b)=(2,2q+1,1,1,1)$ in Lemma~\ref{mine} was obtained by Georghiou and Philippou \cite[p.~35, Eq.~(3.1)]{two}, the case $(n,a,b)=(1,1,1)$ was derived by Borwein et al. \cite{borwein}, and the sums resulting from the case $n=1$ were derived by Flajolet and Salvy \cite{flaj}, and Olaikhan \cite{olaikhan}.

\begin{remark}\label{r1}
	When $p=b=1$, the first term in Lemma~\ref{mine} becomes $\frac{1}{2}\bigl(1+{(-1)}^{q}\bigr)\Li_q(a)\Li_1(1)$,
	which contains the singularity $\Li_1(1)$. Nevertheless, the case $p=b=1$ remains valid because this term is canceled by the endpoint term arising from the sum $\sum_{j=0}^{\lfloor q/2 \rfloor}$:
	\begin{equation*}
		\sum_{j=0}^{\lfloor \frac{q}{2} \rfloor} f(2j)
		=
		\sum_{j=0}^{\lfloor \frac{q-1}{2} \rfloor} f(2j)
		+
		\frac{1}{2}\bigl(1+{(-1)}^{q}\bigr)f(q).
	\end{equation*}
	This cancellation removes the apparent singularity and allows us to formally set $\Li_1(1)=0$. Another point that requires attention is that the $j=0$ terms in the sums $\sum_{j=0}^{\lfloor q/2 \rfloor}$ and $\sum_{j=0}^{\lfloor p/2 \rfloor}$ produce the factor $\Li_0(c)$. When $c=1$, this becomes $\Li_0(1)$, which is interpreted through the analytic continuation of the Riemann zeta function: $\Li_0(1)=\zeta(0)=-\frac12$.
	Therefore, in order for Lemma~\ref{mine} to work properly in Mathematica, it is important to impose the conventions $\Li_1(1)=0$ and $\Li_0(1)=-\frac12$.
\end{remark}

\section{Main results}

\begin{theorem}\label{thm}
	Let $p,n\in\mathbb{Z}_{\geqslant 1}$, $q\in\mathbb{Z}_{\geqslant 2}$, $a=\pm 1$, and $p+q$ odd. Then	
	\begin{align*}
		\sum_{k=1}^\infty\frac{a^k H_{k/n}^{(p)}}{k^q}=
		&\,{(-1)}^{q}\sum_{r=1}^{p}\frac{a}{n^{q-r}}\binom{p+q-r-1}{p-r}\sum_{j=1}^\infty \frac{{(a^n)}^j\mathbf{H}_{nj,a}^{(r)}}{j^{p+q-r}}\\
		&-{(-1)}^{q}\sum_{r=2}^{p}\frac{1}{n^{q-r}}\binom{p+q-r-1}{p-r}\Li_r(a)\Li_{p+q-r}(a^n)\\
		&-{(-1)}^{q}\sum_{\ell=2}^{q-1}\frac{{(-1)}^{\ell}}{n^{q-\ell}}\displaystyle\binom{p+q-\ell-1}{q-\ell}\Li_\ell(a)\zeta(p+q-\ell)\\
		&+\frac{{(-1)}^{q}}{n^{q-1}}\binom{p+q-2}{q-1}\Li_1(a)\{\zeta(p+q-1)-\Li_{p+q-1}(a^n)\},
	\end{align*}
	where $\Li_p(\cdot)$ denotes the polylogarithm, $\zeta(\cdot)$ is the Riemann zeta function, and the sum $\sum_{j=1}^\infty {(a^n)}^j\mathbf{H}_{nj,a}^{(p)}/j^{q}$ is given in Lemma \ref{mine} with $(a,b)$ replaced by $(a^n,a)$. When $a=1$, the final term is understood in the limiting sense as $a\to1$; equivalently, one may set $\Li_1(1)=0$, since the entire term vanishes in this limit.
\end{theorem}

\begin{proof}
	Employing the infinite-series representation of the generalized harmonic numbers \cite[p.~50, Eq.~(1.127)]{ali},
	\begin{equation*}
		H_n^{(p)}
		=\zeta(p)-\sum_{k=1}^{\infty}\frac{1}{{(k+n)}^p},
		\qquad p>1,\quad
		n\notin\mathbb{Z}_{\leqslant -1},
	\end{equation*}
	we have
	\begin{align*}
		\sum_{k=1}^\infty\frac{a^k H_{k/n}^{(p)}}{k^q}=&\sum_{k=1}^\infty\frac{a^k}{k^q}\left(\zeta(p)-\sum_{j=1}^\infty\frac{1}{{(j+k/n)}^p}\right)\\
		=&\Li_q(a)\zeta(p)-n^p\sum_{k=1}^\infty a^k\sum_{j=1}^\infty \frac{1}{k^q{(k+nj)}^p},\quad p>1.
	\end{align*}
	Recalling Lemma~\ref{partial} and interchanging the order of summation, we obtain
	\begin{align*}
		\sum_{k=1}^\infty\frac{a^k H_{k/n}^{(p)}}{k^q}=&\Li_q(a)\zeta(p)-{(-1)}^{q}\sum_{r=1}^{p}\frac{1}{n^{q-r}}\binom{p+q-r-1}{p-r}\sum_{j=1}^\infty \frac{1}{j^{p+q-r}}\sum_{k=1}^\infty  
		\frac{a^k}{{(k+nj)}^r}\\
		&-{(-1)}^{q}\sum_{\ell=1}^{q}\frac{{(-1)}^{\ell} }{n^{q-\ell}}\binom{p+q-\ell-1}{q-\ell}\sum_{j=1}^\infty \frac{ 1
		}{j^{p+q-\ell}}\sum_{k=1}^\infty 
		\frac{a^k}{k^\ell}.
	\end{align*}
	Applying Lemma \ref{m} yields
	\begin{align*}
		\sum_{k=1}^\infty\frac{a^k H_{k/n}^{(p)}}{k^q}=&\Li_q(a)\zeta(p)+{(-1)}^{q}\sum_{r=1}^{p}\frac{a}{n^{q-r}}\binom{p+q-r-1}{p-r}\sum_{j=1}^\infty \frac{{(a^n)}^j\mathbf{H}_{nj,a}^{(r)}}{j^{p+q-r}}\\
		&-{(-1)}^{q}\sum_{r=1}^{p}\frac{1}{n^{q-r}}\binom{p+q-r-1}{p-r}\Li_r(a)\sum_{j=1}^\infty \frac{{(a^n)}^j}{j^{p+q-r}}\\
		&-{(-1)}^{q}\sum_{\ell=1}^{q}\frac{{(-1)}^{\ell}}{n^{q-\ell}}\displaystyle\binom{p+q-\ell-1}{q-\ell}\sum_{j=1}^\infty \frac{1
		}{j^{p+q-\ell}}\sum_{k=1}^\infty 
		\frac{a^k}{k^\ell}.
	\end{align*}
	Finally, using the series representation of the polylogarithm,
	we arrive at
	\begin{align*}
		\sum_{k=1}^\infty\frac{a^k H_{k/n}^{(p)}}{k^q}=&\Li_q(a)\zeta(p)+{(-1)}^{q}\sum_{r=1}^{p}\frac{a}{n^{q-r}}\binom{p+q-r-1}{p-r}\sum_{j=1}^\infty \frac{{(a^n)}^j\mathbf{H}_{nj,a}^{(r)}}{j^{p+q-r}}\\
		&-{(-1)}^{q}\sum_{r=1}^{p}\frac{1}{n^{q-r}}\binom{p+q-r-1}{p-r}\Li_r(a)\Li_{p+q-r}(a^n)\\
		&-{(-1)}^{q}\sum_{\ell=1}^{q}\frac{{(-1)}^{\ell}}{n^{q-\ell}}\displaystyle\binom{p+q-\ell-1}{q-\ell}\zeta(p+q-\ell)\Li_\ell(a).
	\end{align*}
	Observe that the term $\Li_q(a)\zeta(p)$ is canceled by the $\ell=q$ term allowing the case $p=1$. The proof is then completed by isolating the terms $r=1$ and $\ell=1$.
\end{proof}
The case $p=1$ in Theorem~\ref{thm} was recently investigated by Sofo \cite{sofo}. 

\begin{corollary}\label{cor1}
	When $a=1$, Theorem \ref{thm} simplifies to 
\end{corollary}
\begin{align*}
	\sum_{k=1}^\infty\frac{H_{k/n}^{(p)}}{k^q}
	=&\,
	{(-1)}^{q}
	\sum_{r=1}^{p}
	\frac{1}{n^{q-r}}
	\binom{p+q-r-1}{p-r}
	\sum_{k=1}^{\infty}
	\frac{H_{n k}^{(r)}}{k^{p+q-r}}\notag
	\\
	&
	-{(-1)}^{q}
	\sum_{r=2}^{p}
	\frac{1}{n^{q-r}}
	\binom{p+q-r-1}{p-r}
	\zeta(r)\zeta(p+q-r)\notag
	\\
	&
	-{(-1)}^{q}
	\sum_{\ell=2}^{q-1}
	\frac{{(-1)}^{\ell}}{n^{q-\ell}}
	\binom{p+q-\ell-1}{q-\ell}
	\zeta(\ell)\zeta(p+q-\ell),
\end{align*}
which is identical to Sofo's identity \eqref{sofo}.
\begin{corollary}\label{cor2}
	When $a=-1$, Theorem \ref{thm} simplifies to 
	\begin{align*}
		\sum_{k=1}^\infty\frac{{(-1)}^k H_{k/n}^{(p)}}{k^q}=
		&-{(-1)}^{q}\sum_{r=1}^{p}\frac{1}{n^{q-r}}\binom{p+q-r-1}{p-r}\sum_{j=1}^\infty \frac{{({(-1)}^n)}^j\overline{H}_{nj}^{(r)}}{j^{p+q-r}}\\
		&+{(-1)}^{q}\sum_{r=2}^{p}\frac{1}{n^{q-r}}\binom{p+q-r-1}{p-r}\eta(r)\Li_{p+q-r}({(-1)}^n)\\
		&+{(-1)}^{q}\sum_{\ell=2}^{q-1}\frac{{(-1)}^{\ell}}{n^{q-\ell}}\displaystyle\binom{p+q-\ell-1}{q-\ell}\eta(\ell)\zeta(p+q-\ell)\\
		&-\frac{{(-1)}^{q}}{n^{q-1}}\binom{p+q-2}{q-1}\ln(2)\{\zeta(p+q-1)-\Li_{p+q-1}({(-1)}^n)\},
	\end{align*}
	where $\eta(s)=(1-2^{1-s})\zeta(s)$ is the Dirichlet eta function.
\end{corollary}

When $n$ is even, Corollary~\ref{cor2} becomes
\begin{align*}
	\sum_{k=1}^\infty\frac{{(-1)}^k H_{k/n}^{(p)}}{k^q}
	=&-{(-1)}^{q}
	\sum_{r=1}^{p}
	\frac{1}{n^{q-r}}
	\binom{p+q-r-1}{p-r}
	\sum_{j=1}^\infty
	\frac{\overline{H}_{nj}^{(r)}}{j^{p+q-r}}
	\\
	&+{(-1)}^{q}
	\sum_{r=2}^{p}
	\frac{1}{n^{q-r}}
	\binom{p+q-r-1}{p-r}
	\eta(r)\zeta(p+q-r)
	\\
	&+{(-1)}^{q}
	\sum_{\ell=2}^{q-1}
	\frac{{(-1)}^{\ell}}{n^{q-\ell}}
	\binom{p+q-\ell-1}{q-\ell}
	\eta(\ell)\zeta(p+q-\ell).
\end{align*}
Moreover, the case when $n$ is even can also be derived from Corollary~\ref{cor1} by applying the identity
\begin{equation*}
	\sum_{k=1}^\infty {(-1)}^k f(k)
	=
	2\sum_{k=1}^\infty f(2k)
	-
	\sum_{k=1}^\infty f(k),
\end{equation*}
with
$
f(k)=H_{k/2n}^{(p)}/k^q.
$ This gives
\begin{equation*}
	\sum_{k=1}^\infty\frac{{(-1)}^k H_{k/2n}^{(p)}}{k^q}
	=
	\frac{1}{2^{q-1}}
	\sum_{k=1}^\infty\frac{H_{k/n}^{(p)}}{k^q}
	-
	\sum_{k=1}^\infty\frac{H_{k/2n}^{(p)}}{k^q}.
\end{equation*}

When $n$ is odd, Corollary~\ref{cor2} becomes
\begin{align*}
	\sum_{k=1}^\infty\frac{{(-1)}^k H_{k/n}^{(p)}}{k^q}
	=&-{(-1)}^{q}
	\sum_{r=1}^{p}
	\frac{1}{n^{q-r}}
	\binom{p+q-r-1}{p-r}
	\sum_{j=1}^\infty
	\frac{{(-1)}^{j}\overline{H}_{nj}^{(r)}}{j^{p+q-r}}
	\\
	&-{(-1)}^{q}
	\sum_{r=2}^{p}
	\frac{1}{n^{q-r}}
	\binom{p+q-r-1}{p-r}
	\eta(r)\eta(p+q-r)
	\\
	&+{(-1)}^{q}
	\sum_{\ell=2}^{q-1}
	\frac{{(-1)}^{\ell}}{n^{q-\ell}}
	\binom{p+q-\ell-1}{q-\ell}
	\eta(\ell)\zeta(p+q-\ell)
	\\
	&-\frac{2{(-1)}^{q}}{n^{q-1}}
	\binom{p+q-2}{q-1}
	\ln(2)\lambda(p+q-1),
\end{align*}
where
$
\lambda(s)=(1-2^{-s})\zeta(s)
$
is the Dirichlet lambda function.

\begin{theorem}\label{thm2}
	Let $p,q\in\mathbb{Z}_{\geqslant 2}$, $n\in\mathbb{Z}_{\geqslant 1}$, $a=\pm 1$, and $p+q$ odd. Then	
	\begin{align*}
		\int_0^1\frac{\ln^{p-1}(x)\Li_q(ax)}{x(1-x^n)}\,\diff x=&
		-(p-1)!\sum_{r=1}^{p}\frac{a}{n^{p+q-r}}\binom{p+q-r-1}{p-r}\sum_{j=1}^\infty \frac{{(a^n)}^j\mathbf{H}_{nj,a}^{(r)}}{j^{p+q-r}}\\		
		&+(p-1)!\sum_{r=2}^{p}\frac{1}{n^{p+q-r}}\binom{p+q-r-1}{p-r}\Li_r(a)\Li_{p+q-r}(a^n)\\
		&+(p-1)!\sum_{\ell=2}^{q-1}\frac{{(-1)}^{\ell}}{n^{p+q-\ell}}\displaystyle\binom{p+q-\ell-1}{q-\ell}\Li_\ell(a)\zeta(p+q-\ell)\\
		&	-\frac{(p-1)!}{n^{p+q-1}}\binom{p+q-2}{q-1}\Li_1(a)\{\zeta(p+q-1)-\Li_{p+q-1}(a^n)\}\\
		&-\frac{{(-1)}^{p}(p-1)!}{n^p}\{n^p\Li_{p+q}(a)+\Li_q(a)\zeta(p)\},
	\end{align*}
	where $\Li_p(\cdot)$ denotes the polylogarithm, $\zeta(\cdot)$ is the Riemann zeta function, and the sum $\sum_{j=1}^\infty {(a^n)}^j\mathbf{H}_{nj,a}^{(p)}/j^{q}$ is given in Lemma \ref{mine} with $(a,b)$ replaced by $(a^n,a)$.
\end{theorem}

\begin{proof} Using the integral form of the harmonic number \cite[p.~49, Eq.~(1.125)]{ali}
	\begin{equation*}
		H_n^{(p)}=	\zeta(p)-\frac{{(-1)}^{p-1}}{(p-1)!}\int_0^1\frac{x^{n}\ln^{p-1}(x)}{1-x}\,\diff x,\quad p\in\mathbb{Z}_{\geqslant 2},\quad n>-1,
	\end{equation*}
	we have
	\begin{align*}
		\sum_{k=1}^\infty\frac{a^k H_{k/n}^{(p)}}{k^q}=&\sum_{k=1}^\infty\frac{a^k}{k^q}\left(\zeta(p)-\frac{{(-1)}^{p-1}}{(p-1)!}\int_0^1\frac{x^{k/n}\ln^{p-1}(x)}{1-x}\,\diff x\right)\\
		=&\Li_q(a)\zeta(p)-\frac{{(-1)}^{p-1}}{{(p-1)!}}\int_0^1\frac{\ln^{p-1}(x)}{1-x}\sum_{k=1}^\infty\frac{(ax^{1/n})^k}{k^q}\,\diff x\\
		=&\Li_q(a)\zeta(p)-\frac{{(-1)}^{p-1}}{{(p-1)!}}\int_0^1\frac{\ln^{p-1}(x)\Li_q(ax^{1/n})}{1-x}\,\diff x.
	\end{align*}
	Making the change of variable $x^{1/n}\to x$, we get
	\begin{align*}
		\int_0^1\frac{\ln^{p-1}(x)\Li_q(ax^{1/n})}{1-x}\,\diff x=&n^p\int_0^1\frac{x^{n-1}\ln^{p-1}(x)\Li_q(ax)}{1-x^n}\,\diff x\\
		=&n^p\int_0^1\frac{\ln^{p-1}(x)\Li_q(ax)}{x(1-x^n)}\,\diff x-n^p\int_0^1\frac{\ln^{p-1}(x)\Li_q(ax)}{x}\,\diff x.
	\end{align*}
	For the second integral, expanding the polylogarithm into series produces
	\begin{align*}
		\int_0^1\frac{\ln^{p-1}(x)\Li_q(ax)}{x}\,\diff x&=\int_0^1\frac{\ln^{p-1}(x)}{x}\sum_{j=1}^\infty\frac{(ax)^j}{j^q}\,\diff x\\
		&=\sum_{j=1}^\infty\frac{a^j}{j^q}\int_0^1 x^{j-1}\ln^{p-1}(x)\,\diff x\\
		&=\sum_{j=1}^\infty\frac{a^j}{j^q}\frac{{(-1)}^{p-1}(p-1)!}{j^p}\\
		&={(-1)}^{p-1}(p-1)!\Li_{p+q}(a).
	\end{align*}
	Collecting all pieces, we arrive at
	\begin{align*}
		\int_0^1\frac{\ln^{p-1}(x)\Li_q(ax)}{x(1-x^n)}\,\diff x=&\,{(-1)}^{p-1}(p-1)!\Li_{p+q}(a)\\
		&+\frac{{(-1)}^{p-1}(p-1)!}{n^p}\left(\Li_q(a)\zeta(p)-\sum_{k=1}^\infty\frac{a^k H_{k/n}^{(p)}}{k^q}\right).
	\end{align*}
	The sum $\sum_{k=1}^\infty a^k H_{k/n}^{(p)}/k^q$ is given in Theorem~\ref{thm} and the proof is then completed.

\end{proof}


\section{Applications}\label{ex}	
To obtain results from all formulas in this work without encountering divergence issues, it's important to use the following \textit{Mathematica} code:
\begin{verbatim}
	Unprotect[PolyLog];
	PolyLog[0, 1] = -1/2;
	PolyLog[1, 1] = 0;
	Protect[PolyLog];	
	theta[1, 1, s_, r_] := Pi Cot[Pi s/r];
	theta[c_, q_, s_, r_] := LerchPhi[c, q, s/r] +
	c {(-1)}^{q} LerchPhi[c, q, (r - s)/r];	
	formulaExpression // FullSimplify
\end{verbatim}

\subsection{Applications of Corollary~\ref{cor1}}
\begin{example} Let $(p,q,n)=(1,2,2)$. Then
	\begin{equation*}
		\sum_{k=1}^\infty\frac{H_{k/2}}{k^2}=	\frac{11 \zeta (3)}{8}.
	\end{equation*}
\end{example}

\begin{example} Let $(p,q,n)=(2,3,3)$. Then
	\begin{equation*}
		\sum_{k=1}^\infty\frac{H_{k/3}^{(2)}}{k^3}=\frac{233 \zeta (5)}{54}-\frac{17 \pi ^2 \zeta (3)}{54}.
	\end{equation*}
\end{example}

\begin{example} Let $(p,q,n)=(3,4,3)$. Then
	\begin{equation*}
		\sum_{k=1}^\infty\frac{H_{k/3}^{(3)}}{k^4}=-\frac{2 \pi ^3 \zeta (4,\frac{1}{3})}{243 \sqrt{3}}+\frac{2 \pi ^3 \zeta (4,\frac{2}{3})}{243 \sqrt{3}}-\frac{5 \pi ^2 \zeta (5)}{27}+\frac{1111 \zeta (7)}{81}.
	\end{equation*}
\end{example}

\begin{example} Let $(p,q,n)=(5,2,5)$. Then
	\begin{align*}
		\sum_{k=1}^\infty\frac{H_{k/5}^{(5)}}{k^2}=&-\frac{\zeta (2,\frac{1}{5}) \zeta (5,\frac{1}{5})}{50} +\frac{\zeta (2,\frac{4}{5}) \zeta (5,\frac{1}{5})}{50} -\frac{\zeta (2,\frac{2}{5}) \zeta (5,\frac{2}{5})}{50}
		+\frac{\zeta (2,\frac{3}{5}) \zeta (5,\frac{2}{5})}{50}\\
		&+\frac{\zeta (2,\frac{2}{5}) \zeta (5,\frac{3}{5})}{50} -\frac{\zeta (2,\frac{3}{5}) \zeta (5,\frac{3}{5})}{50} +\frac{\zeta (2,\frac{1}{5}) \zeta (5,\frac{4}{5})}{50}
	    -\frac{\zeta (2,\frac{4}{5}) \zeta (5,\frac{4}{5})}{50}\\
		&-\frac{5 \pi ^4 \zeta (3)}{9} -\frac{2 \pi ^2 \zeta (5)}{3}+\frac{39073 \zeta (7)}{25}.
	\end{align*}
\end{example}

\subsection{Applications of Corollary~\ref{cor2}}

\begin{example} Let $(p,q,n)=(1,2,2)$. Then
	\begin{equation*}
		\sum_{k=1}^\infty\frac{{(-1)}^kH_{k/2}}{k^2}=-\frac{3 \zeta (3)}{8}.
	\end{equation*}
\end{example}

\begin{example} Let $(p,q,n)=(2,3,3)$. Then
	\begin{equation*}
		\sum_{k=1}^\infty\frac{{(-1)}^k H_{k/3}^{(2)}}{k^3}=\frac{47 \pi ^2 \zeta (3)}{144}-\frac{3619 \zeta (5)}{864}.
	\end{equation*}
\end{example}

\begin{example} Let $(p,q,n)=(3,4,3)$. Then
	\begin{align*}
		\sum_{k=1}^\infty\frac{{(-1)}^kH_{k/3}^{(3)}}{k^4}=\,\frac{\pi ^3 \zeta (4,\frac{1}{3})}{1944 \sqrt{3}}-\frac{\pi ^3 \zeta (4,\frac{2}{3})}{1944 \sqrt{3}}+\frac{\pi ^3 \zeta (4,\frac{1}{6})}{1944 \sqrt{3}}-\frac{\pi ^3 \zeta (4,\frac{5}{6})}{1944 \sqrt{3}}
		+\frac{\pi ^2 \zeta (5)}{48}-\frac{46027 \zeta (7)}{3456}.
	\end{align*}
\end{example}

\begin{example} Let $(p,q,n)=(5,2,5)$. Then
	\begin{align*}
		\sum_{k=1}^\infty\frac{{(-1)}^kH_{k/5}^{(5)}}{k^2}=&\,\frac{\zeta (2,\frac{2}{5}) \zeta (5,\frac{1}{5})}{200} -\frac{\zeta (2,\frac{3}{5}) \zeta (5,\frac{1}{5})}{200} +\frac{\zeta (2,\frac{1}{10}) \zeta (5,\frac{1}{5})}{200} -\frac{\zeta (2,\frac{9}{10}) \zeta (5,\frac{1}{5})}{200}\\
		&-\frac{\zeta \left(2,\frac{1}{5}\right) \zeta (5,\frac{2}{5})}{200} +\frac{\zeta (2,\frac{4}{5}) \zeta (5,\frac{2}{5})}{200} -\frac{\zeta (2,\frac{3}{10}) \zeta (5,\frac{2}{5})}{200} +\frac{\zeta (2,\frac{7}{10}) \zeta (5,\frac{2}{5})}{200}\\
		&+\frac{\zeta \left(2,\frac{1}{5}\right) \zeta (5,\frac{3}{5})}{200} -\frac{\zeta (2,\frac{4}{5})\zeta (5,\frac{3}{5})}{200}  +\frac{\zeta (2,\frac{3}{10})\zeta (5,\frac{3}{5})}{200}  -\frac{\zeta (2,\frac{7}{10}) \zeta (5,\frac{3}{5})}{200} \\
		&-\frac{\zeta (2,\frac{2}{5}) \zeta (5,\frac{4}{5})}{200} +\frac{\zeta (2,\frac{3}{5}) \zeta (5,\frac{4}{5})}{200} -\frac{\zeta (2,\frac{1}{10}) \zeta (5,\frac{4}{5})}{200} +\frac{\zeta (2,\frac{9}{10}) \zeta (5,\frac{4}{5})}{200} \\
		&-\frac{35 \pi ^4 \zeta (3)}{96}-\frac{5 \pi ^2 \zeta (5)}{16}-\frac{4921293 \zeta (7)}{3200}.
	\end{align*}
\end{example}

\bibliographystyle{amsplain}

\begin{thebibliography}{999}
	
	
\bibitem{andrews}
{\sc G.~E.~Andrews, R.~Askey, and R.~Roy}, \emph{Special Functions},
Cambridge University Press, Cambridge, 1999.

\bibitem{borwein}
{\sc D.~Borwein, J.~M.~Borwein, and R.~Girgensohn}, \emph{Explicit evaluation of Euler sums},
Proc.~Edinburgh Math.~Soc. \textbf{38} (1995), 277--294.
\url{https://carmamaths.org/resources/db90/pdfs/db90-107.00.pdf}

\bibitem{flaj}
{\sc P.~Flajolet and B.~Salvy}, \emph{Euler sums and contour integral representations},
Exp.~Math. \textbf{7} (1998), 15--35.
\url{https://algo.inria.fr/flajolet/Publications/FlSa98.pdf}

\bibitem{two}
{\sc C.~Georghiou and A.~N.~Philippou}, \emph{Harmonic sums and the zeta function},
The Fibonacci Quart. \textbf{21} (1983), 29--36.
\url{https://www.researchgate.net/publication/264239612}

\bibitem{lewin}
{\sc L.~Lewin}, \emph{Polylogarithms and Associated Functions},
North-Holland, New York, 1981.

\bibitem{ali}
{\sc A.~Olaikhan}, \emph{An Introduction to the Harmonic Series and Logarithmic Integrals},
2nd ed., Self-published, Phoenix, AZ, 2023.
\url{https://www.researchgate.net/publication/373488370}

\bibitem{olaikhan2026}
{\sc A.~S.~Olaikhan}, \emph{Closed-form expressions for polylogarithmic integrals and related harmonic sums},
Acta Comment.~Univ.~Tartu.~Math. \textbf{30} (2026), no.~1, 155--177.
\url{https://doi.org/10.12697/ACUTM.2026.30.09}

\bibitem{olaikhan}
{\sc A.~Olaikhan}, \emph{Euler-type sums containing four parameters},
Mathematica Pannonica (published online ahead of print, 2025), Article~314.2025.00016.
\url{https://doi.org/10.1556/314.2025.00016}



\bibitem{hurwitz}
{\sc F.~W.~J.~Olver, D.~W.~Lozier, R.~F.~Boisvert, and C.~W.~Clark},
\emph{NIST Handbook of Mathematical Functions},
Cambridge University Press, New York, 2010.
\url{https://dlmf.nist.gov/}

\bibitem{sofo}
{\sc A.~Sofo}, \emph{Multiple Argument Euler Sum Identities},
Mathematics \textbf{13} (2025), no.~5, Article~839.
\url{https://doi.org/10.3390/math13050839}

\bibitem{sofo2019}
{\sc A.~Sofo}, \emph{General order Euler sums with rational argument},
Integral Transforms Spec.~Funct. \textbf{30}~(12) (2019), 978--991.
\url{https://doi.org/10.1080/10652469.2019.1643851}


\bibitem{choibook}
{\sc H.~M.~Srivastava and J.~Choi}, \emph{Zeta and (q)-Zeta Functions and Associated Series and Integrals},
Elsevier, London, 2012.

\bibitem{xu}
{\sc C.~Xu, Y.~Yan, and Z.~Shi}, \emph{Euler sums and integrals of polylogarithm functions},
J.~Number Theory \textbf{165} (2016), 84--108.
\url{https://doi.org/10.1016/j.jnt.2016.01.025}





\end{thebibliography}

\end{document}